\documentclass[journal abbreviation, manuscript]{copernicus}

\nolinenumbers

\theoremstyle{plain}

\usepackage{overpic}
\theoremstyle{definition}

\usepackage[justification=centering]{caption}
\usepackage{array,caption,threeparttable}
\theoremstyle{remark}

\newcommand{\rP}{\mathrm{P\space }} 
\newcommand{\rE}{\mathrm{E\space }} 

\usepackage{comment}

\usepackage{mathtools}

\DeclarePairedDelimiter\floor{\lfloor}{\rfloor}
\begin{document}

\title{Optimal scheduling of the next preventive maintenance activity for a wind farm}


\Author{Quanjiang}{Yu}
\Author{Michael}{Patriksson}
\Author{Serik}{Sagitov}

\affil{Department of Mathematical Sciences, Chalmers University of Technology and
University of Gothenburg, SE-42 196 Gothenburg, Sweden}




\correspondence{Quanjiang Yu (yuqu@chalmers.se)}

\runningtitle{TEXT}

\runningauthor{TEXT}

\received{}
\pubdiscuss{} 
\revised{}
\accepted{}
\published{}


\firstpage{1}

\maketitle

\begin{abstract}

\noindent
A large part of the operational cost for a wind power farm is due to the cost of  equipment maintenance, especially for offshore wind farms. How to reduce the maintenance cost, and hence increase profitability, is this article’s focus.
It presents a binary linear optimization model whose solution may suggest the wind turbine owners which components, and when, should undergo the next preventive maintenance (PM) replacements. The suggested short-term scheduling strategy takes into account eventual failure events of the multi-component system, in that after the failed system is repaired, the previously scheduled PM plan should be updated, assuming that the restored components are as good as new. 

The optimization algorithm of this paper, NextPM, is tested through numerical case studies applied to a four component model of a wind turbine. The first study addresses the important case of a single component system, used for parameter calibration purposes. The second study analyses the case of seasonal variations of mobilization costs, as compared to the constant mobilization cost setting. Among other things, this analysis reveals a  $35\%$ cost reduction achieved by the NextPM model, as compared to the pure corrective maintenance (CM) strategy. 
The third case study compares the NextPM model with another optimization model - the preventive maintenance scheduling problem with interval costs (PMSPIC), which was the major source of inspiration for this article. This comparison demonstrates that the NextPM model is accurate and much faster in terms of computational time. 

\end{abstract}


\introduction  

Wind energy is one of the lowest-priced renewable energy technologies available today; see \cite{lazard2020lazard}.
A large part of the total cost associated with wind turbines is due to operation and maintenance, amounting to $34\%$ for the fixed-bottom offshore wind turbines, according to \cite{stehly20202019}. To reduce the maintenance cost, one can improve the design of
the components, making them more reliable. One can also reduce the maintenance costs by means of an improved scheduling of the maintenance activities for still functioning components depending on their current age. The latter task is the main concern of this paper, which  proposes an optimization model for preventive maintenance (PM) scheduling of a wind turbine 
or even a farm of wind turbines. 
Notice, that in this paper, by PM activities we don’t mean the practice of regular inspection of the component's condition. Our concern is the optimal planning of preventive \textit{replacements} of the components based on their current age. 

Typically, a maintenance model distinguishes between a corrective maintenance (CM) event,  when a component should be attended after it breaks down, and a PM event, when one or several older components are renewed before they break down, see  the recent survey \cite{lee2016new}.
 An optimal PM scheduling is aimed at reducing the lost production due to the down-time caused by CM events.
 
 There is a multitude of papers devoted to the optimal PM scheduling for multi-component systems, see \cite{werbinska2019technical}.
 The article \cite{jafari2018joint} proposes a joint optimization of the maintenance policy
    and the inspection interval for a multi-unit series system with economic dependence. 
It develops an algorithm aiming at a maintenance policy
    for a multi-component system minimizing the maintenance cost, under the assumption that  one unit of the system is
    subject to condition monitoring, while for
    the other units only the age information is available. 
\cite{tian2014condition} develop a method to quantify the uncertainty of the remaining life length resulting in an effective condition-based maintenance approach to optimal scheduling. 
 
 The article \cite{sarker2016minimizing} looks at opportunistic maintenance (OM), which is a special kind of a PM activity occurring at the time of a CM replacement: replacing still functioning components together with the broken one, may save some mobilization costs. OM activities are shown to be extremely beneficial for the offshore wind farms, due to the large mobilization costs.
 
  In \cite{moghaddam2011sensitivity}, optimization models are developed to determine the optimal PM schedules in repairable and maintainable systems. They show that if mobilization costs are the same irrespective of the number of components to be attended, then multiple simultaneous PM activities become cost-effective. However, their optimization models are nonlinear and non-convex, which makes them computationally hard to solve, see Section 1.3 in \cite{andreasson2020introduction}.
 
The Preventive Maintenance Scheduling Problem with Interval Costs (PMSPIC) model from 
\cite{gustavsson2014preventive} was the major inspiration for this work. The main feature of the PMSPIC model  is the idea of interval cost: given a time interval between two consecutive PM activities, the expected maintenance cost should take into account eventual breakdowns of components during this time interval.
The PMSPIC model has a long computational time, which motivated us to build a new optimization model for PM scheduling of a wing turbine. 


In this paper, we build on the state of the art with a new algorithm, NextPM. Given the current ages of the key components of the system, NextPM computes the best time to perform the next maintenance activity and determines which components should be replaced at that time. The algorithm can be solved in one second, and thus has a potential for being used as a key module in a maintenance scheduling app for wind turbines.

\begin{figure}[h]
\centering
\begin{overpic}[scale=0.55]{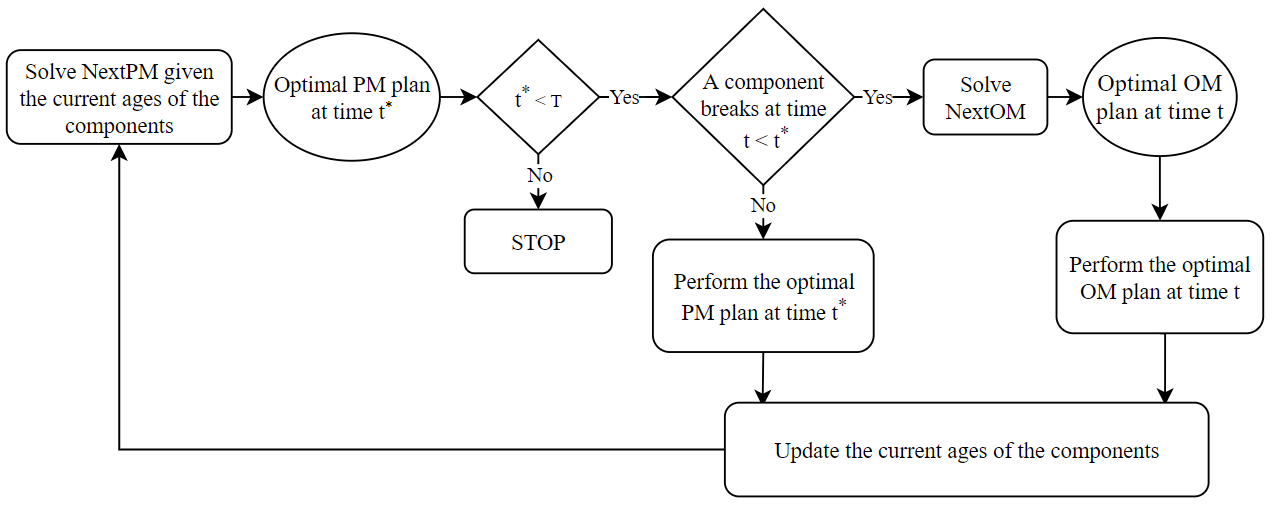}  
\end{overpic}
        \caption{Flow diagram of the optimization algorithm involving NextPM as a major step.}
            \label{flo}
\end{figure}

The paper is organized as follows. Section \ref{alg} presents a novel optimization model for maintenance scheduling of a multi-component system. In the context of wind farm maintenance, each wind turbine is viewed here as a system comprising multiple components such as the gearbox, power generator, rotor, and main bearing. Whenever one of the components is broken, the whole system stops functioning. After the broken component is replaced by a new one, the system resumes its function.
It is assumed that at time $0$ all components of the system are new and that the total lifespan of the system is $T$ units of time. The model has a discrete time setting $t=0,1,\ldots,T,$
where the unit of time can be either a day, or a month, or a year, depending on a concrete application, see \cite{browell2016forecasting} for a maintenance scheduling with only one day ahead. 
In the same Section \ref{alg},  the main result of the paper is summarised as Algorithm 1 aiming at an optimal PM schedule for the time period $[s,T]$ with an arbitrary starting time  $s\in[0,T-1]$. Figure \ref{flo} gives a non-technical description of the algorithm.

The key ingredient of Algorithm 1, the NextPM optimization model, is carefully described in Section \ref{next}. 
Section \ref{num} contains several numerical studies that demonstrate the flexibility of our approach, its accuracy and computational effectiveness. Finally, Section \ref{colu} presents the main conclusions of the paper.


\section{Optimal rescheduling algorithm}\label{alg} 
Consider a system composed of $n$ components characterised by different life length distributions. For the component $j$, it is assumed that its total life length $L_j$ is a random variable having a Weibull distribution with parameters $(\alpha_j,\beta_j)$, so that the corresponding survival function is
\begin{equation}
\rP(L_j>t)=e^{-(\frac{t}{\alpha_j})^{\beta_j}},\quad t\ge0,\quad j=1,\ldots, n;    \end{equation}
see  \cite{guo2009reliability} concerning the use of the Weibull distribution for the modelling of multi-component systems.
The means and variances of the component  life lengths are the following functions of the Weibull parameters
\begin{equation}
    \mu_j=\alpha_j\Gamma(1+\tfrac{1}{\beta_j}),\quad \sigma_j^2=\alpha_j^2\Gamma(1+\tfrac{2}{\beta_j})-\mu_j^2,\quad j=1,\ldots, n.
\end{equation}
Besides the Weibull parameters
$(\alpha_j,\beta_j)$,  $j=1,\ldots,n$, our optimization model 
requires  the following parameters associated with various maintenance costs:
\begin{quote}
$d_t$, the time-dependent mobilization cost for either a PM or CM activity,  $t=0,\ldots,T$,\\
$b_j$, the CM cost of  the component $j=1,\ldots,n$,\\
$c_j$, the PM cost of the component $j=1,\ldots,n$.
\end{quote}
The full set of the model parameters 
$\{d_1,\ldots,d_T,\ (\alpha_1,\beta_1,b_1,c_1),\ldots,(\alpha_n,\beta_n,b_n,c_n),\ \lambda\}$
includes an extra parameter $\lambda$ introduced in Section \ref{cnt} by formula \eqref{G}. 

Suppose that the multi-component system is observed at some time $s\in[0,T-1]$, and the latest maintenance times of components $j=1,\ldots,n$ are known to be $t_j\in[0,s]$, so that at the time $s$, the $n$  components have the effective ages $(s-t_1,\ldots,s-t_n)$. 
The NextPM optimization model described in Section \ref{next} has the input
$(t_1,\ldots,t_n,s,r)$, where $r\in[s+1,T]$ is the end of the current planning period. The output of NextPM is a PM plan specifying the optimal time $\tau\in [s+1,r+1]$ of the next PM event, as well as the set of components $\mathcal{P}\subset\{1,\ldots,n\}$ which should be maintained at the time $\tau$. In particular, the output $\tau=r+1$ means that no PM activity should be scheduled during the planning period $[s+1,r]$ implying that the set $\mathcal{P}$ is empty.

The NextPM model is the key module of the following Algorithm 1 aiming at the long-term PM scheduling until the end-time $T$, at which the whole system is expected to be dismantled, see Figure \ref{flo} for a flowchart illustrating the major steps of Algorithm 1.
\begin{algorithm}[H]
\caption{Optimal rescheduling algorithm}
\hspace{0cm}Input  $t_1,\ldots,t_n,s,r$

\hspace{0.5cm}Start

\hspace{1cm}Solve NextPM\{$t_1,\ldots,t_n,s,r$\}

\hspace{1cm}Output $\tau,\ \mathcal{P}$, where $\mathcal{P}\subset\{1,\ldots,n\}$ is the set of components subject to PM activities at time $\tau$

\hspace{1cm}If $\tau<T$

\hspace{2cm}If $\text{a failure during the period } (s,\tau] \text{ damages component } i \text{ at time } u_{i}$

\hspace{3cm}Set $u:=\floor{u_{i}}$

\hspace{3cm}Solve NextOM\{$i,t_1,\ldots,t_n,u$\}

\hspace{3cm}Output $\mathcal{O}\subset\{1,\ldots,n\}$ is the set of components subject to OM activities at time $u+1$

\hspace{3cm}Perform CM of component $i$ at time $u+1$

\hspace{3cm}Perform PM of each component $j\in\mathcal{O}$ at time $u+1$

\hspace{3cm}Update $r:=\min(u+1+r-s,T)$, $s:=u+1$

\hspace{3cm}Update $t_j:=u+1,\ j\in\mathcal{O}\cup\{i\}$

\hspace{2cm}Else
\hspace{0.3cm}Perform PM of each component $j\in\mathcal{P}$ at time $\tau$

\hspace{3cm}Update $r:=\min(\tau+r-s,T)$, $s:=\tau$, $t_j:=s,\ j\in\mathcal{P}$

\hspace{2cm}End

\hspace{2cm}Go to Start

\hspace{0.5cm}Stop

\end{algorithm}
Algorithm 1 relies on a rescheduling procedure, where each NextPM step covering $r-s$ units of the planning time is accompanied by a NextOM module. The latter is a modification of the NextPM step, see Section \ref{OM}, which addresses 
the possibility of a component failure before the planned PM, followed by an OM activity.

\section{An optimal plan for  the next preventive maintenance}\label{next}

\begin{figure}[h]
\centering
\begin{overpic}[scale=0.6]{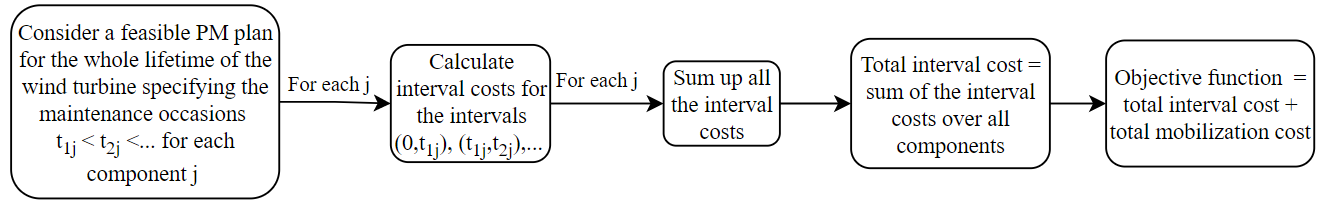}  
\end{overpic}
        \caption{A flow diagram demonstrating how the PMSPIC calculates the objective function for a given feasible maintenance plan.}
            \label{flow}
\end{figure}

\begin{figure}[h]
\centering
\begin{overpic}[scale=0.55]{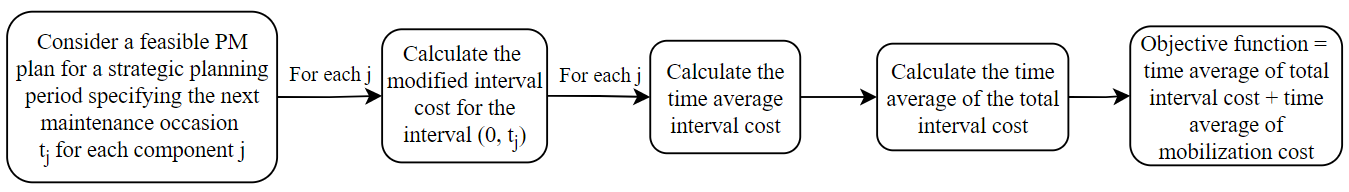}  
\end{overpic}
        \caption{A flow diagram demonstrating how NextPM model calculates the objective function for a given feasible maintenance plan with $s=0$.}
            \label{flon}
\end{figure}

This section sets up the optimization model NextPM, which is the key ingredient of Algorithm 1 summarised in Section \ref{alg}. The optimization model PMSPIC of \cite{gustavsson2014preventive} was a major motivation for NextPM, and we start by comparing these two approaches using Figures \ref{flow} and  \ref{flon} which illustrate two different definitions of the objective functions for two optimization models in question. 
The main difference between PMSPIC and NextPM model is that while PMSPIC generates a maintenance plan for the whole lifetime of the wind turbine, the NextPM model produces an optimal schedule only for the next PM activity. To this end, PMSPIC looks into the total maintenance cost, while NextPM aims at minimizing the time average maintenance cost.

\subsection{NextPM model}\label{nPM}
The purpose of the NextPM model is to produce an optimal PM plan for the period $[s+1,r]$, where the planning timespan $r-s$ is chosen so that it is reasonable to expect at most one PM event during time $r-s$.
For a given planning period $[s+1,r]\subset[0,T]$, an {\it $(s,r)$-plan} is defined as a collection $({\boldsymbol{z}},{\boldsymbol{x}}^ 1,\ldots,{\boldsymbol{x}}^ n)$ of vectors
$${\boldsymbol{z}}=(z_{s+1},\ldots, z_{r+1}),\qquad {\boldsymbol{x}}^ j=(x_{s+1}^ j,\ldots, x_{r+1}^ j),\quad j=1,\ldots,n,$$ 
with binary coordinates
$z_{t},\ x_{t}^ j\in \{0,1\},$
which satisfy the following linear conditions:
\begin{align}
 &\sum_{t=s+1}^{r+1}x_t^ j= 1,\quad j=1,\ldots,n,\\
&x_{t}^ j\le z_t,\quad t=s+1,\ldots, r+1,\ j=1,\ldots,n.
\end{align}
For $t=s+1,\ldots r$, the equality $x_t^ j=1$ means that
\begin{quote}
according to the $(s,r)$-plan, component $j$ should undergo a PM replacement at time $t$, provided no component failure during the time period $[s+1,t]$. 
\end{quote}
In contrast, the equality $x_{r+1}^ j=1$  means that
\begin{quote}
according to the $(s,r)$-plan, no PM activity should involve component $j$ during the time period $[s+1,r]$. 
\end{quote}
On the whole system level, the equality $z_t=1$ means that 
\begin{quote}
according to the $(s,r)$-plan, at least one component should undergo a PM replacement at time $t$, provided no component failure during the time period $[s+1,t]$.
\end{quote}
and the equality $z_{r+1}=1$ means that 
\begin{quote}
according to the $(s,r)$-plan, no PM activity is scheduled for the time period $[s+1,r]$. 
\end{quote}

The NextPM optimization model is built around the objective function 
\begin{equation}
f({\boldsymbol{z}},{\boldsymbol{x}}^ 1,\ldots,{\boldsymbol{x}}^ n)=\sum_{t=s+1}^{r+1} \frac{1}{t-s}\Big(d_tz_t+ c_{s,t}^ 1x^ 1_{t}+\ldots+c_{s,t}^ nx^ n_{t}\Big),
 \label{f}
\end{equation}
where $d_tz_t$ stands for the mobilization cost and the terms $c_{s, t}^ j$ are the so called interval costs defined in Section \ref{cnt}. 
The objective function 
\eqref{f} can be viewed as the {\it time-average} maintenance cost per time unit according to the $(s,t)$-plan $({\boldsymbol{z}},{\boldsymbol{x}}^ 1,\ldots,{\boldsymbol{x}}^ n)$. 

Let $(\bar{\boldsymbol{z}},\bar{\boldsymbol{x}})$ be the solution to the linear optimization problem aimed to
\begin{equation}
 \text{minimize} \quad f({\boldsymbol{z}},{\boldsymbol{x}}^ 1,\ldots,{\boldsymbol{x}}^ n),
\end{equation}
over all $(s,t)$-plans subject to the linear constraints
\begin{equation}
D_{s,t}^ jx^ j_{t}\geq 0,\quad t=s+1,\ldots ,r,\ j=1,\ldots,n,
 \label{D}
\end{equation}
where
$D_{s,t}^ j$ is defined in Section \ref{Dnt} as the PM benefit  for the component $j$ at time $t$.
Then the NextPM algorithm $(\tau, \mathcal N)$ computes the optimal time of the next PM by
\[\tau=  \min_{j} \{\arg\max_t \bar x_t^ j\},\]
and determines the set of the components that should undergo the maintenance activities at time $\tau$ using
\begin{align*}
 \mathcal N&=\left\{
\begin{array}{ll}
\{j: \bar x_{\tau}^ j=1,\ j=1,\ldots,n\}  &  \text { if } \tau\le r,  \\
 \emptyset&   \text { if } \tau= r+1.
\end{array}
\right.
\end{align*}

\subsection{Definition of modified interval costs $c_{s,t}^ j$}\label{cnt}

Here we deal with the term $c_{s,t}^ j$ appearing in the the objective function \eqref{f} of the optimization model NextPM.
The main idea is to define  $c_{s,t}^ j$ as the fixed PM cost $c_j$ plus the expected additional costs due to eventual failures of the component $j$ occurring prior to the planned PM activity at time $t$. 

To this end, consider $n$ independent sequences of renewal times with a delay by letting $U_{s,0}^ j=s$, 
\begin{equation}
    U_{s,1}^ j=t_j+L_{1j},\quad L_{1j}\stackrel{d}{=}\{L_j|L_j>s-t_j\},
\end{equation}
where $\stackrel{d}{=}$ means equality in distribution (conditional distribution in the above formula), and 
\begin{equation}
    U_{s,i+1}^ j=U_{s,i}^ j+L_{ij},\quad L_{ij}\stackrel{d}{=}L_j,\quad \text{for } i=2,3,\ldots, 
\end{equation}
assuming that the random variables $(L_{ij})$ are mutually independent. Notice that in the important particular case $s=0$, this definition simplifies, so that for each $j$, the sequence $\{U_{0,i}^ j\}_{i\ge0}$ describes a renewal process without a delay.

Treating $U_{s,1}^ j, U_{s,2}^ j, \ldots$ as the sequence of consecutive failure times of the component $j$, put
\begin{equation}
\begin{aligned}
    c_{s,t}^ j:=c_j+\rE\Bigg(\sum_{i=1}^{\infty}1_{\{U_{s,i}^ j\le t\}}G_j(U_{s,i-1}^ j,L_{ij},t-s)
    \Bigg),
    \end{aligned}
\label{A}
\end{equation}
where the cost functions
\begin{equation}
G_j(s,u,v)=b_j+d_{s+u}-(\tfrac{u}{v})^{\lambda}\left(c_j+d_{s+v}\right),\quad 0\le u\le v,
\label{G}
\end{equation}
involve a new parameter $\lambda>0$ assumed to be independent of $j=1,\ldots,n$. 
The definition of the cost function  \eqref{G} further develops the key idea of Section $5.1$ in \cite{gustavsson2014preventive}. 
It describes the additional cost implied by an eventual breakdown of component $j$ before the planned PM activity. 

The expression \eqref{G} is suggested as a compromise between two extreme cases: a failure at the start of the planning period, $u=0$, and a failure just before the planned PM replacement, $u=v=t-s$. If $u$ is close to $0$, then the failure at time $s+u$ will not change the PM plan, implying that  the much smaller additional cost 
\begin{equation}
G_j(s,0,t)=b_j+d_{s}  
\label{max}
\end{equation}
is the sum of the CM cost $b_j$ and the mobilization cost $d_s$ at time $s$. On the other hand, if $u$ is close to $v=t-s$, then the additional cost
\begin{equation}
 G_j(s,t,t-s)=b_j-c_j   
\label{min}
\end{equation}
is simply the difference between the CM and PM costs. For $u\in(0,v)$, the expression in the right-hand side of \eqref{G} produces an additional cost which lies between the extreme values \eqref{max} and \eqref{min}. The role of the parameter $\lambda$ is to control to what extent the proximity of the failure time to the planned PM time influences the extra costs. For example, if $\lambda=1$ the intermediate cost is found by a linear extrapolation. 

\subsection{Definition of $D^ j_{s,t}$}\label{Dnt}

The constraint \eqref{D} arises as a check-up step to ensure that a suggested PM at time $t$ brings some  benefit, as compared to a simple strategy when no PM is performed. With the PM-free strategy, the total maintenance cost (including mobilization costs) for the component $j$ during the period $[s,T]$ would be
\begin{equation*}
   \rE\left[\sum_{i=1}^{\infty}1_{\{U_{s,i}^ j\leq T\}}\Big(b_j+d_{U_{s,i}^ j}\Big)\right].  
\end{equation*}
   Alternatively, if the plan  is to perform a PM for  the component $j$ at time $t$, and then to perform replacements of the component $j$ whenever it breaks down, then the total cost would be 
   \begin{equation*}
   c_{s,t}^ j+\rE\left[\sum_{i=1}^{\infty}1_{\{t+U_{0,i}^ j\leq T\}}\Big(b_j+d_{t+U_{0,i}^ j}\Big)\right].    
   \end{equation*}
   Taking into account the difference between these two total costs 
   \begin{equation}
       D^ j_{s,t}=\rE\left[\sum_{i=1}^{\infty}1_{\{U_{s,i}^ j\leq T\}}\Big(b_j+d_{U_{s,i}^ j}\Big)\right]-c_{s,t}^ j-\rE\left[\sum_{i=1}^{\infty}1_{\{t+U_{0,i}^ j\leq T\}}\Big(b_j+d_{t+U_{0,i}^ j}\Big)\right],
   \end{equation}
we conclude that the planned PM of the component $j$ at time $t$ is justified only if 
$D^ j_{s,t}\ge 0$.

\subsection{Complete optimization model of NextPM}
Here we put together the complete optimization model of the NextPM step:
\begin{subequations}
\begin{alignat*}{5}
\mathop{\rm minimize}\ &&f({\boldsymbol{z}},{\boldsymbol{x}}^ 1,\ldots,{\boldsymbol{x}}^ n) &:=\sum_{t=s+1}^{r+1} \frac{1}{t-s}\Big(d_tz_t+ c_{s,t}^ 1x^ 1_{t}+\ldots+c_{s,t}^ nx^ n_{t}\Big), \span\span
        \\
\mathop{\rm subject\ to}& &\sum_{t=s+1}^{r+1}x_t^ j&= 1,\hspace{8mm}\quad j=1,\ldots,n,\\
&&z_t&\geq x_{t}^ j,\hspace{6mm}\quad t=s+1,\ldots, r+1,\ j=1,\ldots,n,\\
&&D_{s,t}^ jx^ j_{t}&\geq 0,\hspace{8mm}\quad t=s+1,\ldots, r,\ j=1,\ldots,n,\\
        & &z_{t}&\in \{0,1\},\quad t=s+1,\ldots, r+1,\\
&& x_{t}^{j}&\in \{0,1\},\hspace{0.5mm}\quad t=s+1,\ldots, r+1,\ j=1,\ldots,n.
\end{alignat*}
\end{subequations}

\subsection{NextOM model}\label{OM}

The NextOM step of Algorithm 1 is a specialised version of the NextPM step described above. The input vector of the NextOM algorithm
$$(i,t_1,\ldots,t_{i-1},t_{i+1},\ldots,t_n,s),$$   
treats $i$ as the label of the component whose failure at some time during $[s,s+1)$ has triggered the OM planning step. For a pair $\{s,i\}$, an {\it $\{s,i\}$-plan} is any set of vectors $({\boldsymbol{z}},{\boldsymbol{x}}^ 1,\ldots,{\boldsymbol{x}}^ n)$ whose components are two-dimensional vectors 
\begin{equation}
    {\boldsymbol{z}}=(z_{s+1}, z_{s+2}),\qquad {\boldsymbol{x}}^ j=(x_{s+1}^ j, x_{s+2}^ j),\quad j=1,\ldots,n,
\end{equation}
with binary coordinates
$z_{t},\ x_{t}^ j\in \{0,1\}$
satisfying the following linear conditions
\begin{align}
 &\sum_{t=s+1}^{s+2}x_t^ j= 1,\quad j=1,
 \ldots,n,\\
&x_{s+1}^{(i)}= 1,\\
&z_t\geq x_{t}^ j,\quad t=s+1,s+2,\quad j=1,\ldots,n.
\end{align}
Observe that necessarily, $z_{s+1}=1$. 

The NextOM optimization model uses the following modified version of the objective function \eqref{f}:
\begin{equation}
f_i({\boldsymbol{z}},{\boldsymbol{x}}^ 1,\ldots,{\boldsymbol{x}}^ n)=\sum_{t=s+1}^{s+2} \frac{1}{t-s}\Big(d_tz_t+ \sum_{j
\neq i}c_{s,t}^ jx^ j_{t}\Big),
\end{equation}
where $c_{s,t}^ j$ is defined in Section \ref{cnt}. 
Let $(\bar{\boldsymbol{z}},\bar{\boldsymbol{x}})$ be the solution to the linear optimization problem to
\begin{equation}
\text{minimise }\quad f_i({\boldsymbol{z}},{\boldsymbol{x}}^ 1,\ldots,{\boldsymbol{x}}^ n)    
\end{equation}
over all $\{s,i\}$-plans subject to the linear constraints
\begin{equation}
D_{s,s+1}^ jx^ j_{s+1}\geq 0,\quad j=1,\ldots,i-1,i+1,\ldots,n,
\end{equation}
where
$D_{s,t}^ j$ is defined in Section \ref{Dnt}.
The output of the NextOM  is given by the set 
\begin{align*}
 \mathcal O&=\{j: \bar x_{\tau}^ j=1,\ j=1,,\ldots,i-1,i+1,\ldots,n\},
\end{align*}
consisting of the labels of the components which will be opportunistically maintained along with the component $i$ undergoing a CM activity.

\section{Numerical studies} \label{num}
The three case studies analysed in this section treat a wind turbine as a system represented by four components.
They are all based on the parameter values taken from the paper \cite{tian2011condition}, see  Table \ref{t1},
\begin{table}[h]
\centering
\caption{Key parameters for a four-component system.}
  \begin{tabular}{|l|c|c|c|c|r|r|}\hline
 \small Component&$j$&CM cost $b_j$ $(\$ 1000)$&PM cost $c_j$ $(\$ 1000)$&$\beta_j$&$\alpha_j$ (months)&$\mu_j$ (months)\\\hline 
    \small Gearbox&1&202&46.75&3&80&71.4\\
         \small Rotor&2&162&36.75&3&100& 89.9\\
              \small Generator&3&150&33.75&2&110&97.5\\
\small Main bearing&4&110&23.75&2&125&110.8\\\hline
\end{tabular}
\label{t1}
\end{table}
where the cost unit is $1000$ USD and the time unit is 1 month. The lifetime of the wind turbine is assumed to be $30$ years. 
This implies the parameter value $T=360$ months. As to other model parameters, it is assumed that 
\begin{description}
\item[ ] $s=0$ which implies that all four components initially are as good as new,
\item[ ] $r=60$  months, see Section \ref{cs1} for motivation,
\item[ ] $\lambda=3$, based on \cite{gustavsson2014preventive}.
\end{description}

Comparing the characteristics of four wind turbine components shown in Table \ref{t1}, it is important to observe a strict ordering from the perspective of the associated PM costs. For example, consider components 1 and 2. The rightmost column says that the expected life length of the gearbox is smaller by 18.5 months, which on its own, suggests a higher rate of replacements for the component 1. But even the other two parameters, CM cost and PM cost, are ordered in a way $b_1>b_2$, $c_1>c_2$, which is favorable for more frequent replacements of component 1 compared to component 2.

All computational tests are performed on an Intel 2.40 GHz dual
core Windows PC with 16 GB RAM. The mathematical optimization
models are implemented in AMPL IDE (version 3.5); 
the model components \eqref{A} and \eqref{D} are calculated by Matlab (version R2015b), 
and then the optimization problems are solved using
CPLEX (version 12.8).

\subsection{Study 1: focusing on a single component at a time}\label{cs1}
If $n=1$, $d_t\equiv d$, and $s=0$, the objective function \eqref{f} takes the form
\begin{equation}
f({\boldsymbol{x}})=\sum_{t=1}^{r+1} a_tx_{t},\quad a_t=\frac{d+ c_{t}}{t},    
\label{ata}
\end{equation}
where given a sequence of independent random variables $L_i\stackrel{d}{=}L$ with $L$ having a Weibull $(\alpha,\beta)$ distribution,
\begin{equation}
    c_t=c+\rE\left(\sum_{i=1}^\infty 1_{\{L_1+\ldots+L_i\le t\}}\Big[b+d-(\tfrac{L_i}{t})^\lambda(c+d)\Big]\right).
\end{equation}
In the single component setting, coefficient $a_t$ in \eqref{ata} describes the monthly maintenance cost if the next PM is planned at time $t$ (assuming that at time 0 the component was as good as new). In this section, we analyze the behavior of the function $a_t$ under some realistic model parameters. It turns out in the current setting, that minimising the objective function  \eqref{f} is equivalent to
minimising $a_t$ over $t=1,\ldots,r+1$, and moreover, the constraint \eqref{D} can effectively be disregarded.  As a result of this analysis, we propose $r=60$ months as a practical length of the planning period for our algorithm.

\begin{figure}[h]
\centering

\begin{overpic}[scale=0.6]{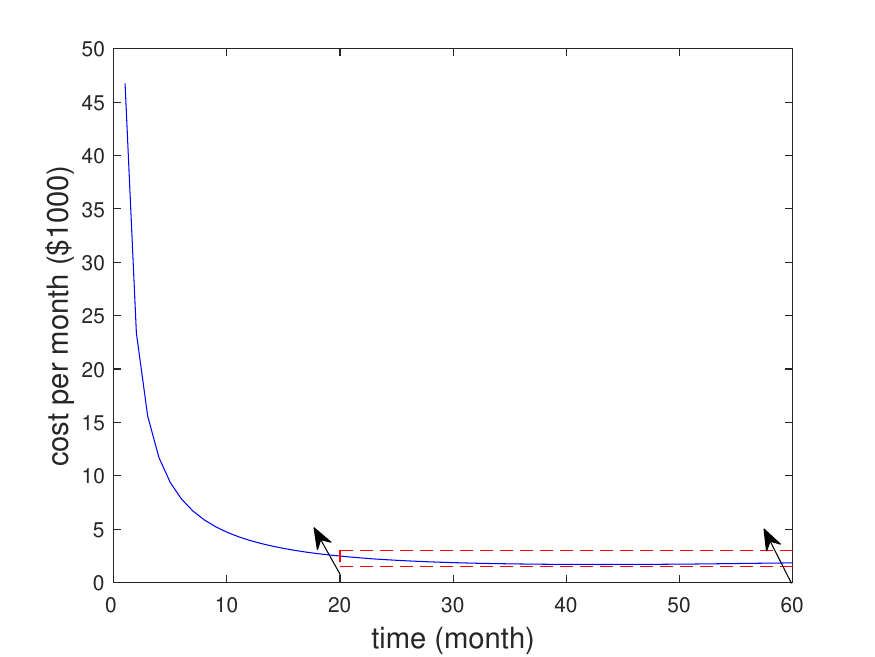}  
\put(23,16){\includegraphics[scale=.4]{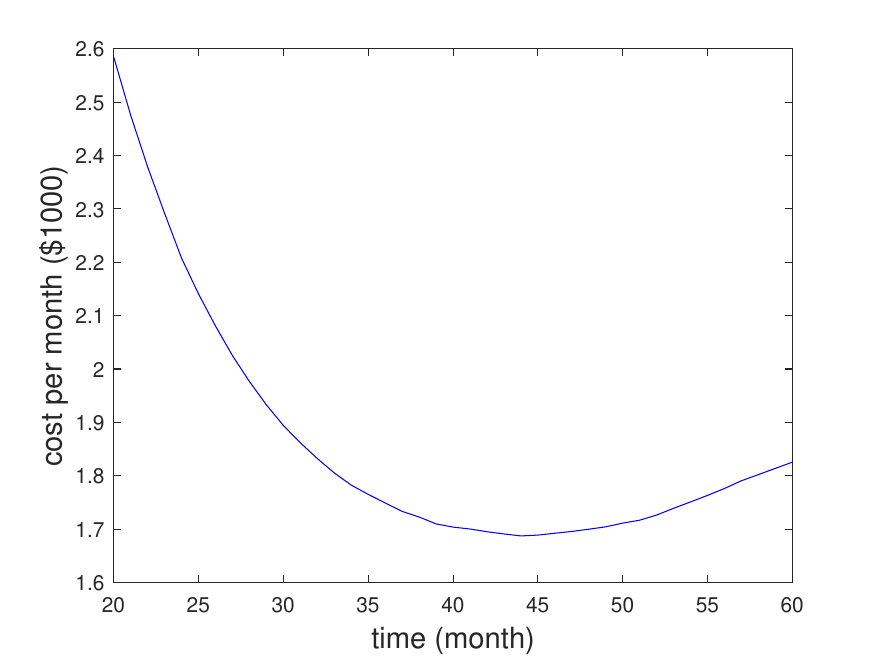}} 
\end{overpic}
        \caption{Monthly maintenance cost for the gearbox with the mobilization cost $d=\$5000$. }
            \label{PM1}
\end{figure}

Figure \ref{PM1} presents a typical profile for the monthly maintenance cost $a_t$ as a function of the time $t$ of the next PM planned activity. The inset of Figure \ref{PM1} clearly shows that the best time for next PM is at $\tau=43$ given the mobilization cost of $d=\$5000$. The maintenance cost in this case is $a_{43}=\$1700$ per month. 

\begin{figure}[h]
\centering
\includegraphics[width=0.5\textwidth]{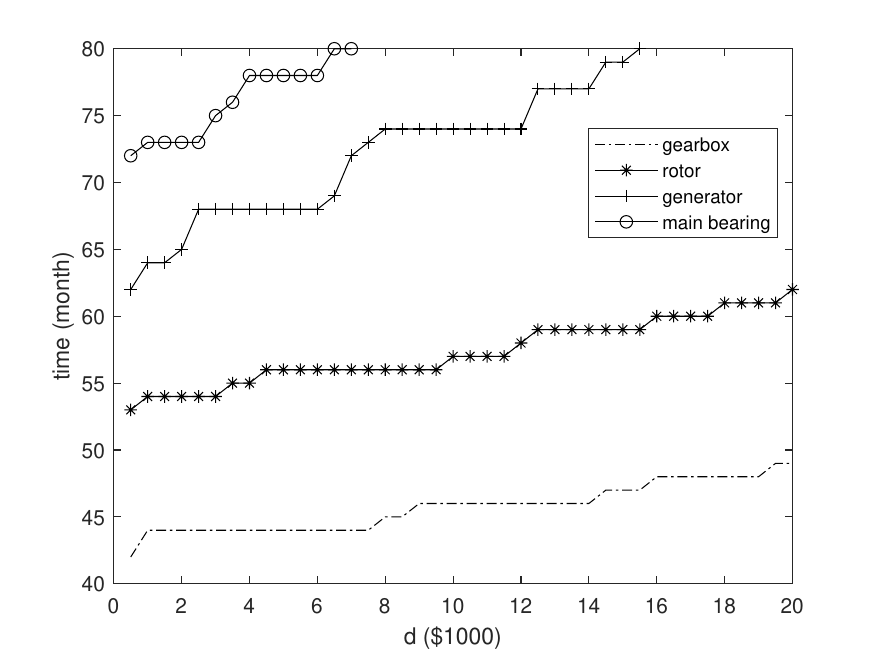}
        \caption{The optimal next PM time $\tau$ as a function of the mobilization cost $d$ for different single-component systems.}
            \label{PM2}
\end{figure}

The same value $\tau=43$ can be also seen on the lowest among four lines depicted on Figure \ref{PM2} if parameter $d$, shown on the horizontal axis, takes value 5. The gearbox line on Figure \ref{PM2} displays larger values of $\tau$ for higher mobilization costs $d$. The same pattern is seen for the other three components taken one at a time.

Now we are ready to explain how the results of our analysis justify the proposed value $r=60$ months for the length of the next PM planning period. The ideal choice of $r$ must satisfy two contadicting requirements. On one hand, $r$ should not be very small to avoid too many NextPM steps in Algorithm 1 advising for no PM activities during the next planning period. On the other hand, a smaller value of the parameter $r$ would significantly reduce the computational time of the NextPM model. As a compromise solution, we choose $r$ in such way that at least one PM activity is expected to be scheduled during the planning horizon. Since gearbox is expected to be replaced most often, referring to Figures \ref{PM1} and \ref{PM2}, we take the value of $r=60$ months for our case studies. 

\subsection{Study 2: seasonal effects}
In this section, we study how different mobilization costs $d_t$ result in different optimal PM schedules. Part A presents a baseline study of a pure CM strategy with no PM activities. Part B deals with seasonally changing $d_t$ around the average value $\bar d =\$10000$. Part C takes up a similar case with a lower average mobilization cost $\bar d = \$5000$. \\

Part A. 

Consider the simplest wind turbine maintenance strategy when the PM option is ignored and a CM activity is performed whenever a turbine component breaks down. This baseline case study will help us to evaluate how much can be saved by introducing PM planning. 

The total cost associated with the pure CM strategy is estimated based on the random number of failures over the time interval [$0,T]$ for all $n$ components 
\begin{equation}
    F(T)=\sum_{j=1}^n\rE\left(\sum_{i=1}^\infty 1_{\{V^ j_i\le T\}}(d_{V^ j_i}+b_j)\right)=\sum_{j=1}^n\int_0^T(d_u+b_j)dH_j(u),
\end{equation}
where $H_j$ are the corresponding renewal functions
\begin{equation}
H_j(t)=\rE\left(\sum_{i=1}^\infty 1_{\{V^ j_i\le t\}}\right),\quad t>0,\ j=1,\ldots,n.    
\end{equation}
According to the standard renewal theory, see for example  \cite{grimmett2020probability}, for large values of $T$,
\begin{equation}
\frac{F(T)}{T}\approx\sum_{j=1}^n\frac{1}{T\mu_j}\int_0^T(d_u+b_j)du=\sum_{j=1}^n\frac{\bar d+b_j}{\mu_j},    
\end{equation}
where
$
\bar d=\frac{d_1+\ldots+d_T}{T}.    
$
Applying this approximation to the four-component model of the wind turbine, the monthly maintenance costs for the pure CM strategy are computed to be $\$7396$  for $\bar d=\$5000$, and $\$7618$ for $\bar d=\$10000$.\\

Part B. 

To address the seasonal effects of the mobilization costs $d_t$, the following mobilization costs (in thousands of USD) for different months in a year are used:
\begin{center}
  \begin{tabular}{cccccccccccc}
Jan&Feb&Mar&Apr&May&Jun&Jul&Aug&Sep&Oct&Nov&Dec\\\hline 
15&13&11&9&7&5&5&7&9&11&13&15
\end{tabular}
\end{center}
It is assumed that the mobilization costs for different months are different, but for the same month at different years are the same. In this case, the average mobilization cost is $\bar d=\$10000$. 
The given monthly costs are based on a discussion with the experts affiliated with the Swedish Wind Power Technology Centre (SWPTC).
Table \ref{t3} summarises the results produced by the NextPM algorithm applied to the following three settings:
\begin{description}
\item[Winter start] scenario: if the wind turbine starts functioning in January, then the mobilization costs $d_t$ (in thousands of USD) follow the following periodical dynamics over time $t=1,2,\ldots$:
$$d_1=15, \ d_2=13,\ldots, \ d_{12}=15, \ d_{13}=15, \ d_{14}=13,\ldots$$
\item[Summer start] scenario: if the wind turbine starts functioning in July, then the mobilization costs $d_t$ are taking the values (in thousands of USD)
 $$d_1=5, \ d_2=7,\ldots, \ d_{12}=5, \ d_{13}=5, \ d_{14}=7,\ldots$$
\item[Constant mobilization cost] scenario has no seasonal effect in that for each month $t$, the mobilization cost $d_t$ is the same: 
 $d_1=10, \ d_2=10, \ d_3=10,\ldots$ (in thousands of USD).
\end{description}
\begin{table}[h]
\centering
\caption{Summary of the NextPM results for $\bar d=\$10000$.}
  \begin{tabular}{|l|c|c|c|c|c|c|}\hline
Component $j$&1&2&3&4&Corresponding month&Monthly maintenance cost\\\hline 
Winter start&54&54&54&54&Jun&\$5010  \\
Summer start&49&49&49&49&Jul&\$4979 \\
Constant mobilization cost &52&52&52&52&--&\$5061  \\\hline
\end{tabular}
\label{t3}
\end{table}
Our results suggest (as a consequence of high mobilization costs)  to perform PM to all four components at a certain time,  irrespective of the scenario. With the summer start setting, the average monthly maintenance cost is somewhat lower. Notice that in all of the seasonal settings, the proposed PM activities are scheduled for summer months (having lower mobilization costs). 
Observe that all three monthly averages $(\$5010,\$4979,\$5061)$ are much lower than the baseline value $\$7618$ obtained in Part A.\\

Part C. 

In this section, the mobilization costs are halved to contrast the results of Part B, so that $\bar d=\$5000$ and $d_t$ take the following values (in thousands of USD) depending on which month of the year lies behind the time parameter $t$:

\begin{center}
\small 
  \begin{tabular}{cccccccccccc}
Jan&Feb&Mar&Apr&May&Jun&Jul&Aug&Sep&Oct&Nov&Dec\\\hline 
7.5&6.5&5.5&4.5&3.5&2.5&2.5&3.5&4.5&5.5&6.5&7.5
\end{tabular}
\end{center}
\normalsize
The new results presented in Table \ref{t2} are drastically different from the  results of Part B.

\begin{table}[h]
\centering
\caption{Summary of the NextPM results for $\bar d=\$5000$.}
  \begin{tabular}{|l|c|c|c|c|c|c|}\hline
Component $j$&1&2&3&4&Corresponding month&Monthly maintenance cost\\\hline 
Winter start&43&x&x&x&Jul&\$4876 \\
Summer start&48&48&x&x&Jul&\$4863 \\
Constant mobilization cost &50&50&50&50&--&\$4964 \\\hline
\end{tabular}
\label{t2}
\end{table}
According to Table \ref{t2}, in the Winter start setting, the optimal next PM plan suggests a PM activity on month $43$ only for the component 1, the gearbox. With the seasonal mobilization cost, the next PM is always planned during the summer since the mobilization cost is low then. Again, the most economic among the three scenarios is to start in the summer time, with the optimal plan being to perform the next PM activity on month $48$ by replacing the components 1 and 2.

Comparison of the results of Parts B and C with those of Part A shows that implementation of the PM planning reduces maintenance costs by 35\%.

\subsection{A performance comparison with PMSPIC }

In this case study, we compare the outputs of the NextPM model and the optimization model PMSPIC. A comparison of the NextPM model with the PMSPIC model is not a straightforward exercise, since the latter produces a maintenance plan for the whole lifespan $[0,T]$ of the multi-component system in question. To make a fair comparison, we characterize both approaches in terms of  the time average  maintenance costs.
The following three tables summarise the results for three values of the constant mobilization cost $d$:

\small 
\begin{table}[h]
\centering
\caption{Outputs of the NextPM and PMSPIC models for $d=\$1000$.}
 \noindent  \begin{tabular}{l|c|c|c|c||c|r|r}
$d=\$1000$&1&2&3&4&Monthly maintenance cost&Matlab&AMPL\\\hline 
NextPM&43&x&x&x&\$4731 &49 sec&0.01 sec\\
PMSPIC &49&x&x&x&\$4746 &135 sec&19.64 sec
\end{tabular}
\end{table}

\begin{table}[h]
\centering
\caption{Outputs of the NextPM and PMSPIC models for $d=\$5000$.}
\begin{tabular}{l|c|c|c|c||c|r|r} 
$d=\$5000$&1&2&3&4&Monthly maintenance cost&Matlab&AMPL\\\hline 
NextPM&50&50&50&50&\$4964 &54 sec&0.01 sec\\
PMSPIC&51&51&51&51&\$4881  &132 sec&51.56 sec 
\end{tabular}
\end{table}

\begin{table}[h]
\centering
\caption{Outputs of the NextPM and PMSPIC models for $d=\$10000$.}
  \begin{tabular}{l|c|c|c|c||c|r|r}
$d=\$10000$&1&2&3&4&Monthly maintenance cost&Matlab&AMPL\\\hline 
NextPM&52&52&52&52&\$5061 &55 sec&0.01 sec\\
PMSPIC&50&50&50&50&\$5037  &134 sec&87.57 sec
\end{tabular}
\end{table}

\normalsize
\noindent Tables 4-6 reveal that the next PM schedules produced by NextPM and PMSPIC are quite similar. The observed small differences in the maintenance costs do not imply that PMSPIC gives better solutions, since NextPM calculates the maintenance costs within a different modelling framework. 
The main advantage of NextPM compared to PMSPIC is in the computational speed. The effectiveness of the algorithms is reported in the two rightmost columns. The “Matlab” column gives the time it takes to generate the main parameters of the model. For the NextPM the number of parameters is much smaller, and they are $c_{s,t}^j$, $D_{s,t}^j$. The “AMPL” column gives the time it takes to solve the optimization model.  For example, if $d= \$10000$, the NextPM optimization runs $10000$ times faster than the PMSPIC optimization.

For $d=\$5000$, the NextPM calculations are performed with the time unit being three days. The results are rather similar to those obtained for  the time unit 1 month. Solving this problem with AMPL has required a time increase from $0.01$ to $0.08$ seconds caused by a ten-fold increase of the number of the time steps. The corresponding increase in the AMPL time  for the PMSPIC model was much more dramatic: it takes more than 11 hours to solve the full optimization problem.

\section{Conclusions}\label{colu}
This article introduces a new NextPM optimization model aiming at PM scheduling for a wind turbine viewed as multi-component system. Which of the components should undergo PM replacements first is decided based on the information on the component ages. Compared to the PMSPIC model from \citet{gustavsson2014preventive} that generates a maintenance plan for the
whole lifetime of the wind turbine, NextPM model produces an optimal schedule only for the next PM activity. By
focusing on a shorter planning horizon and implementing a different model structure, we succeeded in substantial reduction of the computational time.

NextPM is tested with three case studies based on the data for four components of the wind turbine taken from \citet{tian2011condition}. Under the seasonal variation, our results show that PM activities should be always scheduled in the summer time. This is due to the lower mobilization costs during the summer months.  When the NextPM model is compared to the pure CM strategy, it is found that around $35\%$ of the maintenance costs can be saved by applying the NextPM model. 
The third case study compares the performances of NextPM and PMSPIC algorithms 
demonstrating accuracy of the NextPM model despite of being much less complex than PMSPIC. 

In this paper our NextPM model is applied to a system of four components belonging to a single wind turbine. However, we claim that our approach can handle the case of, say, ten turbines with $80$ components in total (the computational time required by our algorithm grows linearly with the increased number of components, while PMSPIC’s computational time grows exponentially fast.)

In the future, we plan to use NextPM as a key module in a maintenance scheduling app for wind power turbines, after adding a module for processing condition monitoring data.

\appendixfigures  

\appendixtables   


\authorcontribution{Quanjiang Yu developed the theoretical formalism, performed the analytic calculations and performed the numerical
simulations. Serik Sagitov, Quanjiang Yu , and Michael Patriksson  contributed to the final version of the manuscript. } 

\competinginterests{There is no significant competing financial, professional, or personal interests that might have influenced the performance or presentation of the work described in this manuscript.} 


\begin{acknowledgements}
\noindent We acknowledge
the financial support from the Swedish Wind Power Technology Centre at Chalmers, and from the Swedish Research Council (Dnr.\ 2014-5138). Special thanks to the director of SWPTC, professor Ola Carlson, for his constructive recommendations. The valuable comments of four reviewers helped considerably in improving the quality of our manuscript.
\end{acknowledgements}

 \bibliography{paper.bib}

\end{document}